\documentclass[twoside,reqno]{lt6}
\usepackage{cite}
\usepackage{amssymb,amsmath}
\usepackage{times}
\def\z{z}
\def\w{w}
\def\u{u}
\def\spr{}
\def\tt{\otimes}                               

\def\wti{\widetilde}

\def\di{\partial}
\def\Di{\partial}

\newcommand\llb{[\![}                          
\newcommand\rrb{]\!]}                          
\newcommand\llbl{(\!(}                          
\newcommand\rrbl{)\!)}                          
\def\vac{|0\rangle}                            

\def\i{{\mathrm{i}}}     
\def\CC{\mathbb{C}}       
\def\ZZ{\mathbb{Z}}       
\def\QQ{\mathbb{Q}}       
\def\NN{\mathbb{N}}       
\def\NNN{{\ZZ_+}}
\def\CCC{{\CC^\times}}
\def\al{\alpha}                         
\def\be{\beta}
\def\ga{\gamma}
\def\Ga{\Gamma}
\def\de{\delta}
\def\ett{\eta}
\def\oom{\omega}
\def\OM{\Omega}
\def\io{\iota}
\def\De{\Delta}
\def\ep{\varepsilon}

\def\Om{\Omega}
\def\ph{\varphi}
\def\si{\sigma}

\def\h{{\mathfrak{h}}}      

\def\I{\mathcal{I}}    
\def\Y{\mathcal{Y}}    
\def\QF{\mathcal{F}}   

\DeclareMathOperator{\Res}{Res}

\DeclareMathOperator{\id}{id}

\DeclareMathOperator{\End}{End}

\DeclareMathOperator{\Ext}{Ext}

\DeclareMathOperator{\coh}{H}

\renewcommand\Re{{\mathrm{Re}}}

\newcommand{\skp}{\medskip\noindent}
\def\qed{\nopagebreak\nolinebreak\quad\nolinebreak$\Box$}
\newcommand{\Proof}{{\textit{Proof.\; }}}
\newcommand{\Theorem}[1]{{\textbf{Theorem{#1}.}}}
\newcommand{\Lemma}[1]{{\textbf{Lemma{#1}.}}}
\newcommand{\Proposition}[1]{{\textbf{Proposition{#1}.}}}
\newcommand{\Corollary}[1]{{\textbf{Corollary{#1}.}}}
\newcommand{\Definition}[1]{{\textbf{Definition{#1}.}}}
\newcommand{\Example}[1]{{\textbf{Example{#1}.}}}
\newcommand{\Examples}[1]{{\textbf{Examples{#1}.}}}
\newcommand{\Remark}[1]{{\textit{Remark{#1}.}}}
\newcommand{\thref}[1]{Theorem \ref{#1}}
\newcommand{\prref}[1]{Proposition \ref{#1}}
\newcommand{\leref}[1]{Lemma \ref{#1}}
\newcommand{\coref}[1]{Corollary \ref{#1}}
\newcommand{\deref}[1]{Definition \ref{#1}}
\newcommand{\exref}[1]{Example \ref{#1}}
\newcommand{\reref}[1]{Remark \ref{#1}}
\newcommand{\seref}[1]{Sect.\ \ref{#1}}
\begin{document}
\sloppy \raggedbottom
\setcounter{page}{1}

\newpage
\setcounter{figure}{0}
\setcounter{equation}{0}
\setcounter{footnote}{0}
\setcounter{table}{0}
\setcounter{section}{0}



\title{Generalized Vertex Algebras}

\runningheads{Bojko Bakalov, Victor G.~Kac}{Generalized Vertex Algebras}

\begin{start}


\coauthor{Bojko Bakalov}{1},
\coauthor{Victor G.~Kac}{2},

\address{Department of Mathematics, North Carolina State University, 
Box 8205, Raleigh, NC 27695, USA; bojko\_bakalov@ncsu.edu}{1}
\address{Department of Mathematics, MIT, Cambridge, MA 02139, USA; 
kac@math.mit.edu}{2}


\begin{Abstract}
We give a short introduction to generalized vertex algebras, 
using the notion of polylocal fields. 
We construct a generalized vertex algebra associated to a vector space
$\h$ with a symmetric bilinear form. 
It contains as subalgebras all lattice vertex algebras of rank
equal to $\dim \h$ and all irreducible representations of these 
vertex algebras.
\end{Abstract}

\end{start}

\section{Introduction}\label{s1}

A \emph{vertex algebra} is essentially the same as a chiral algebra
in two-dimensional conformal field theory \cite{BPZ,G,DMS}.
In mathematics,
vertex algebras arose naturally in the representation theory
of infinite-dimensional Lie algebras and in the construction
of the ``moonshine module'' for the Monster simple finite group
\cite{B,FLM}. 

Some of the most important vertex (super)algebras 
are the vertex (super)algebras
$V_Q$ associated to \emph{integral lattices} $Q$ 
\cite{FK,B,FLM,K,LL}. 
If the lattice $Q$ is not necessarily integral,
one gets on $V_Q$ the structure of a \emph{generalized vertex algebra}
as introduced in \cite{FFR,DL,M} 
(under the name ``vertex operator para-algebra'' in \cite{FFR}).
This notion includes as special cases the notion of a ``vertex
superalgebra'' (called just a ``vertex algebra'' in \cite{K})
and the notion of a ``colored vertex algebra'' from \cite{X}.
The ``parafermion algebras'' of Zamolodchikov and Fateev \cite{ZF1,ZF2}
are also closely related to generalized vertex algebras
(see \cite[Chapter 14]{DL}).
The theory of generalized vertex algebras (and further
generalizations) is developed in detail in the monograph \cite{DL},
and important examples of generalized vertex algebras are
constructed in \cite{FFR,DL,M,GL}.
The treatment of \cite{FFR,DL,M} is centered around 
a ``Jacobi identity,'' which generalizes the Jacobi identity
from \cite{FLM} (the latter is equivalent to the Borcherds identity
from \cite{K}). 

In the present paper we give a short introduction to 
generalized vertex algebras, utilizing the approach of 
\cite{BN} (for $D=1$)
which is based on the notion of \emph{polylocal fields}.
The notion of \emph{locality} plays an important role
in the theory of (generalized) vertex algebras
(see \cite{G,DL,K,Li1,GL,LL}). It is natural both
from a physical and from a mathematical point of view
to extend it to polylocality of fields in several variables.
Our definition of a generalized vertex algebra is 
slightly more general than the ones from \cite{FFR,DL,M} 
in that we allow a more general grading condition,
and we do not assume an action of the Virasoro algebra
(cf.\ \cite{GL}).
In particular, in \cite{DL,M} only \emph{rational lattices} $Q$
give rise to generalized vertex algebras $V_Q$.
In our setting, $Q$ is allowed to be the whole \emph{vector space},
which leads to a new generalized 
vertex algebra $V_\h$ associated to any 
vector space $\h$ with a 
symmetric bilinear form. 
It is remarkable that $V_\h$ contains as subalgebras the
lattice vertex (super)algebras $V_Q$ 
for all integral lattices $Q \subseteq \h$,
as well as all of their irreducible representations.

The paper is organized as follows.
In \seref{s2} we recall the notion of a \emph{quantum field}
and its generalization that corresponds to taking non-integral powers
of the formal variable. We also introduce fields in several variables and
(generalized) \emph{polylocality}, following \cite{BN}.

In the presence of a grading by an abelian group $Q$, 
we define \emph{parafermion fields} in \seref{s3.1}.
Then in \seref{s3.2} we show that a 
translation covariant, local and complete system of parafermion fields
can be extended uniquely to a \emph{state-field correspondence}.
This result implies Uniqueness and Existence Theorems that
generalize those of \cite{G,FKRW,K} (see also \cite{Li1,GL,LL}).
The definition of a \emph{generalized vertex algebra} is
given in \seref{s3.3} in terms of (generalized) locality. 
In \seref{s3.6} we introduce a natural action of the 
\emph{cohomology group} $\coh^2(Q,\CCC)$ on the isomorphism classes of
generalized vertex algebras. We also show how in a special case the notion of
a generalized vertex algebra reduces to that of a \emph{$Q$-graded}
vertex superalgebra.
In \seref{s3.4} we discuss the \emph{operator product expansion}
of local parafermion fields, and we prove \emph{formal associativity}
and \emph{commutativity} relations generalizing those of \cite{FB,LL}.
In \seref{s3.5} we derive the 
(generalized) \emph{Jacobi identity} (= \emph{Borcherds identity}),
thus showing that our 
definition is a generalization of the ones from \cite{FFR,DL,M}.
The exposition of Sections \ref{s2}, \ref{s3.2}, \ref{s3.4} and \ref{s3.5}  
follows closely \cite{BN}.

In \seref{s4} we introduce the notion of a 
\emph{module} over a generalized vertex algebra, following \cite{DL}.
We show that
the notion of a \emph{twisted module} over a vertex (super)algebra
\cite{FFR,D2} (cf.\ \cite{DK}) 
is a special case of the notion of a module
over a generalized vertex algebra, as observed in \cite{Li2}.

In \seref{s5} we construct the generalized vertex algebra $V_\h$
associated to a \emph{vector space} $\h$ with a 
symmetric bilinear form, and we discuss its subalgebras $V_Q$
associated to \emph{integral lattices} $Q\subseteq\h$.

\section{Quantum fields in several variables and polylocality}\label{s2}

\subsection{Spaces of formal series}\label{s2.1}

We first fix some notation to be used throughout the paper.
Let $\NN=\{1,2,\dots\}$ and $\NNN=\{0,1,2,\dots\}$.
By $\z,\w$, etc., we will denote formal commuting variables.
All vector spaces are over the field $\CC$ of complex numbers.
We will denote by $V$ a vector space, and by 
$V[\z]$ (respectively, $V \llb \z \rrb$)
the space of polynomials (respectively, formal power series)
in $\z$ with coefficients in~$V$.

We will identify the subsets $\Ga \subseteq \CC/\ZZ$ with subsets
$\Ga \subseteq \CC$ that are \emph{$\ZZ$-invariant}, i.e., that satisfy
$\Ga+\ZZ \subseteq \Ga$. Given such a subset $\Ga$, we denote by
$V \llb\z\rrb\z^\Ga$ the space of all finite sums of the form
$\sum_i \psi_i(\z) \z^{d_i}$, where $d_i \in \Ga$ and
$\psi_i(\z) \in V \llb \z \rrb$.
Another way to write the elements of $V \llb\z\rrb\z^\Ga$
is as infinite sums
$\sum_n f_n \z^{n}$, where $f_n \in V$ and $n$
runs over the union of finitely many sets of the form
$\{ d_i + \NNN \}$ with $d_i \in \Ga$.
In particular, $V \llb\z\rrb\z^\ZZ$ is exactly the space
$V \llbl\z\rrbl \equiv V \llb\z\rrb [\z^{-1}]$ 
of formal Laurent series.

We denote by $V \llb\z, \z^\Ga\rrb$ the space of all 
formal infinite series 
$\sum_{n\in\Ga} f_n \z^{n}$ with $f_n \in V$.
For $\Ga=\ZZ$ this coincides with the space $V \llb\z, \z^{-1} \rrb$ 
of formal power series in $\z$, $\z^{-1}$.
Note that $V \llb\z\rrb\z^\Ga$
is a $\CC\llbl\z\rrbl$-module, while $V \llb\z, \z^\Ga\rrb$ 
is only a $\CC[\z,\z^{-1}]$-module, $V \llb\z\rrb\z^\Ga$ being
its submodule.
In addition, both $V \llb\z\rrb\z^\Ga$
and $V \llb\z, \z^\Ga\rrb$ are equipped with the usual
action of the derivative~$\di_\z$.

In the same way, we introduce spaces of formal series 
in several variables
$V \llb \z_1,\dots,\z_s \rrb \z_1^{\Ga_{\! 1}} 
\dotsm \z_s^{\Ga_{\! s}}$
and 
$V \llb \z_1,\z_1^{\Ga_{\! 1}}, \dots,\z_s,\z_s^{\Ga_{\! s}} \rrb$.
The latter consists of all infinite sums of the form
\begin{equation*}
\sum_{n_j \in \Ga_{\! j}} 
f_{n_1,\dots,n_s} \z_1^{n_1} \dotsm \z_s^{n_s} \,,
\qquad f_{n_1,\dots,n_s} \in V \,,
\end{equation*}
while the former is its subspace consisting of such sums with
each $n_j$ running over the union of finitely many sets of the form
$\{ d_{i,j} + \NNN \}$ with $d_{i,j} \in \Ga_j$.
Let us point out that 
$V \llb \z_1,\z_2 \rrb \z_1^{\Ga_{\! 1}} \z_2^{\Ga_{\! 2}}$
consists of series in which the powers of $\z_1$ and $\z_2$
are uniformly bounded from below. In contrast, elements of the space
$\bigl( V \llb \z_1 \rrb \z_1^{\Ga_{\! 1}} \bigr)
\llb \z_2 \rrb \z_2^{\Ga_{\! 2}}$
have powers of $\z_2$ that are
bounded from below but the powers of $\z_1$ are possibly unbounded.

For $N \in\CC$, we define the formal expansions
\begin{align}
& \label{e2.3}
\begin{aligned}
\io_{\z_1,\z_2} (\z_1-\z_2 &)^N := 
e^{-\z_2 \spr \di_{\z_1}} \z_1^N
\\
&= \sum_{j\in\NNN} \binom{N}{j} \z_1^{N-j} (-\z_2)^j
\in \bigl( \CC\llb\z_1\rrb \z_1^{\Ga} \bigr) \llb\z_2\rrb \,,
\end{aligned}
\\ 
& \label{e2.4}
\begin{aligned}
\io_{\z_2,\z_1} (\z_1-\z_2 &)^N := 
e^{\pi\i N} e^{-\z_1 \spr \di_{\z_2}} \z_2^N
\\
&= e^{\pi\i N} \sum_{j\in\NNN} \binom{N}{j} (-\z_1)^j \z_2^{N-j}
\in \bigl( \CC\llb\z_2\rrb \z_2^{\Ga} \bigr) \llb\z_1\rrb \,,
\end{aligned}
\end{align}
where $\Ga=N+\ZZ$.
Note that the spaces in the right-hand sides of Eqs.\ \eqref{e2.3}
and \eqref{e2.4} are modules over the localized ring
$\CC\llb\z_1,\z_2\rrb_{\z_1,\z_2}$. More generally, it makes sense
to multiply Eq.\ \eqref{e2.3} by an element of 
$V \llb\z_1,\z_2\rrb \z_1^{\Ga_{\! 1}} \z_2^{\Ga_{\! 2}}$
thus producing an element of 
$\bigl( V \llb\z_1\rrb \z_1^{N+\Ga_{\! 1}} \bigr) 
\llb\z_2\rrb  \z_2^{\Ga_{\! 2}}$,
and similarly for Eq.\ \eqref{e2.4}.
In this way, we extend the above expansions by linearity to maps
\begin{align*}
\io_{\z_1,\z_2} \colon
V \llb\z_1,\z_2\rrb \z_1^{\Ga_{\! 1}} \z_2^{\Ga_{\! 2}} \z_{12}^{\, N}
& \, \to \,
\bigl( V \llb\z_1\rrb \z_1^{N+\Ga_{\! 1}} \bigr) 
\llb\z_2\rrb  \z_2^{\Ga_{\! 2}} \,,
\\
\io_{\z_2,\z_1} \colon
V \llb\z_1,\z_2\rrb \z_1^{\Ga_{\! 1}} \z_2^{\Ga_{\! 2}} \z_{12}^{\, N}
& \, \to \,
\bigl( V \llb\z_2\rrb \z_2^{N+\Ga_{\! 2}} \bigr) 
\llb\z_1\rrb  \z_1^{\Ga_{\! 1}} \,,
\end{align*}
where $\z_{12} := \z_1-\z_2$.

Obviously, when $N \in\NNN$, expansions \eqref{e2.3} and \eqref{e2.4}
are equal to each other and coincide with the binomial expansion of 
$(\z_1-\z_2)^N$. Similarly, we define expansions of
$(\z_1+\z_2)^N$ for $N\in\CC$.
One has the following analog of \emph{Taylor's formula}:
\begin{equation}\label{e2.15}
\io_{\z_1,\z_2} \, \io_{\z_{12},\z_2} \, f(\z_{12}+\z_2,\z_2)
= \io_{\z_1,\z_2} \, f(\z_1,\z_2) \,,
\quad \z_{12} := \z_1-\z_2 \,,
\end{equation}
for every localized series
$f(\z_1,\z_2) = g(\z_1,\z_2) \, \z_{12}^{\, N}$
with $N\in\CC$ and
$g(\z_1,\z_2) \in \CC\llb\z_1,\z_2\rrb \z_1^\CC \z_2^\CC$.
Indeed, it is enough to prove Eq.\ \eqref{e2.15} for
$f(\z_1,\z_2) = \z_1^M$ with $M\in\CC$,
in which case it follows from \eqref{e2.3}
(see e.g.\ Proposition 2.2 from \cite{BN} for more details).

\subsection{Polylocal quantum fields}\label{s2.2}

We define a \emph{quantum field}
in $m$ variables $\z_1,\dots,\z_m$ (or just an \emph{$m$-field} for short)
to be a linear map from $V$ to the space
$V \llb \z_1,\dots,\z_m \rrb \z_1^\CC \dotsm \z_m^\CC$.
Alternatively, an $m$-field $A (\z_1,\dots,\z_m)$ can be viewed
as a formal series from 
$(\End V) \llb \z_1,\z_1^\CC, \dots,\z_m,\z_m^\CC \rrb$
with the property that for every $v \in V$ one has:
\begin{equation*}
A (\z_1,\dots,\z_m) v = \sum_i \,
\z_1^{d_{i,1}} \dotsm \z_m^{d_{i,m}} 
\, \psi_i (\z_1,\dots,\z_m) \,,
\end{equation*}
where the sum is finite, $d_{i,j} \in \CC$
and $\psi_i \in V \llb \z_1,\dots,\z_m \rrb$.

If $A$ is an $m$-field, then for every partition
\begin{equation*}
\{1,\dots,m\} =
J_1 \sqcup \dotsm \sqcup J_r
\qquad \text{(disjoint union)} , 
\end{equation*}
the \emph{restriction}
\begin{equation}\label{e2.7}
\wti A (\u_1,\dots,\u_r) v := 
\bigl( A (\z_1,\dots,\z_m) v \bigr) 
\big|_{ \z_j := u_s \,\;\text{for}\,\; j \in J_s }
\end{equation}
makes sense and defines an $r$-field.

We will assume that $V$ is endowed with
an endomorphism $T$ (called \emph{translation operator})
and with a vector $\vac$ (called \emph{vacuum vector}), such that
$T \vac = 0$. An $m$-field $A$
is called \emph{translation covariant\/} if
\begin{equation*}
T A (\z_1,\dots,\z_m) - A (\z_1,\dots,\z_m) \, T
= \sum_{k=1}^m \, \di_{z_k} A (\z_1,\dots,\z_m) \,.
\end{equation*}

Let us point out that a product
$A (\z_1,\dots,\z_m) B (\z_{m+1},\dots,\z_{m+n})$
of two fields is \emph{not} a field in general.
Indeed, by the above definition, for every $v \in V$ we have
\begin{align*}
A (\z_1,\dots,\z_m) v & \, \in
V \llb \z_1,\dots,\z_m \rrb 
\z_1^\CC \dotsm \z_m^\CC
\,, 
\\ 
B (\z_{m+1},\dots,\z_{m+n}) v & \, \in
V \llb \z_{m+1},\dots,\z_{m+n} \rrb 
\z_{m+1}^\CC \dotsm \z_{m+n}^\CC
\,,
\end{align*}
which implies that
$A (\z_1,\dots,\z_m) B (\z_{m+1},\dots,\z_{m+n}) v$
belongs to the space
\begin{equation*}
\bigl(
V \llb \z_1,\dots,\z_m \rrb 
\z_1^\CC \dotsm \z_m^\CC
\bigr)
\llb \z_{m+1},\dots,\z_{m+n} \rrb 
\z_{m+1}^\CC \dotsm \z_{m+n}^\CC
\,.
\end{equation*}
In general, elements of the latter space have unbounded powers of 
$\z_1,\dots,\z_m$. As a consequence, the restriction of the above product
for coinciding arguments is not well defined in general.
We will show below that one can ``regularize'' this product
to make a field if the following definition is satisfied 
(see~\cite{BN}).

\skp
\Definition{ \ref{s2.2}}
An $m$-field $A$ and an $n$-field $B$ are called
\emph{mutually local\/} if there exist complex numbers 
$N_{i,j}$ $(i=1,\dots,m$, $j=1,\dots,n)$
and $\ett\in\CCC$ such that
\begin{align}\label{e2.12}
\Bigl(
\prod_{i=1}^m & \prod_{j=m+1}^{m+n} \!
\io_{\z_i,\z_j} (\z_i-\z_j)^{N_{i,j}} 
\Bigr)
A (\z_1,\dots,\z_m) B (\z_{m+1},\dots,\z_{m+n})
\\ \notag
=& \, \ett \,
\Bigl(
\prod_{i=1}^m \prod_{j=m+1}^{m+n} \!
\io_{\z_j,\z_i} (\z_i-\z_j)^{N_{i,j}} 
\Bigr)
B (\z_{m+1},\dots,\z_{m+n}) A (\z_1,\dots,\z_m) 
\,.
\end{align}

A $1$-field that is local with respect to itself is usually called
a \emph{local} field;
a $2$-field that is local with respect to itself is called
a \emph{bilocal} field.
An $m$-field, for general $m$, which is local with respect to itself,
is called a \emph{polylocal} field.

The following important result
is a straightforward variation of Theorem 4.1 from~\cite{BN}.

\skp
\Proposition{ \ref{s2.2}}
{\it{
Let\/ $A (\z_1,\dots,\z_m)$ and\/ $B (\z_{m+1},\dots,\z_{m+n})$
be an\/ $m$-field and an\/ $n$-field, 
respectively, which are mutually local as above.

\medskip

{\rm (a)}
Every restriction of
$A$ 
for coinciding arguments is also a field and
is mutually local with respect to $B$ 
$($see Eq.~\eqref{e2.7}$)$.

\medskip

{\rm (b)}
If the field $A$ 
is translation covariant, then
its restrictions for coinciding arguments 
are also translation covariant fields.

\medskip

{\rm (c)}
If $A$ 
is translation covariant, then
$A (\z_1,\dots,\z_m) \vac \in
V \llb \z_1,\dots,\z_m \rrb$.

\medskip

{\rm (d)}
Every partial derivative $\di_{z_k} A$
is a field and is mutually local with respect to $B$. 
If the field $A$ is translation covariant, then
$\di_{z_k} A$ is also translation covariant.

\medskip

{\rm (e)}
The left-hand side $F_{A,B}$ of\/ Eq.~\eqref{e2.12}
is an $(m+n)$-field. If the fields $A$ 
and $B$ 
are local with respect to a $p$-field $C$,
then $F_{A,B}$ is also local with respect to $C$.
If both fields $A$ and $B$ are translation covariant, then
$F_{A,B}$ is also translation covariant.
}}

\skp
\Remark{ \ref{s2.2}}
The above Proposition remains valid if one considers 
a more general notion of polylocality,
where in Eq.\ \eqref{e2.12} the product $\prod (\z_i-\z_j)^{N_{i,j}}$
is replaced with some function of the differences $\z_i-\z_j$.
One can replace that product with an even more general function
if the translation covariance condition is modified accordingly.
See also \cite{Li3}, where related ideas $($in the case $m=n=1)$
were applied in the investigation of ``quantum vertex algebras.''

\subsection{Operator product expansion}\label{s2.3}

As a corollary of \prref{s2.2},
every $m$-field $A (\z_1,\dots,\z_m)$ can be expanded
in $1$-fields as follows (see \cite{BN}).
Consider for $v \in V$ the formal expansion
\begin{align}
\notag
\iota_{\z,\w_1} & \dotsm \, \iota_{\z,\w_{m-1}} \,
A (\z+\w_1,\dots,\z+\w_{m-1},\z) v
\\
\label{e2.13}
:= & \,
\exp( \w_1 \spr \Di_{\z_1} + \dotsm + \w_{m-1} \spr \Di_{\z_{m-1}} )
\, A (\z_1,\dots,\z_m) v
\big|_{ \z_1 = \dotsm = \z_m = \z  }
\\ 
\notag
\in & \,
\bigl( V \llb \z \rrb \z^\CC \bigr)
\llb \w_1,\dots,\w_{m-1} \rrb \,.
\end{align}
This is a formal power series in $\w_1,\dots,\w_{m-1}$ with
coefficients of the form
$\psi_i (\z) v \in V \llb \z \rrb \z^\CC$
for some uniquely defined fields $\psi_i (\z)$
($i$ running over some index set).
All $\psi_i (\z)$ are fields because
they are obtained from $A (\z_1,\dots,\z_m)$ by the operations
of differentiation and restriction.
If, in addition,
$A$ is translation covariant and is local with respect to
some other fields $B$, $C$, etc., then all the fields $\psi_i (\z)$
are also translation covariant and local with respect to $B$, $C$, etc.

Formal expansion \eqref{e2.13} is called
the \emph{operator expansion} of $A (\z_1,\dots,\z_m)$.
Applying this expansion to the field $F_{A,B}$
given by the left-hand side of Eq.~\eqref{e2.12},
we get what is called the \emph{operator product expansion} (OPE) of two
mutually local fields $A$ and $B$.
The following simple observation will be useful in the sequel.

\skp
\Remark{ \ref{s2.3}}
It follows from \prref{s2.2}(c) that for $v=\vac$
the right-hand side of Eq.\ \eqref{e2.13} is just
the Taylor series expansion of
\begin{equation*}
A (\z+\w_1,\dots,\z+\w_{m-1},\z) \vac
\in
V \llb \z, \w_1,\dots,\w_{m-1} \rrb \,.
\end{equation*}
This implies that the linear span of all coefficients of
$A (\z_1,\dots,\z_m) \vac$ coincides with the linear span
of all coefficients of all $\psi_i (\z) \vac$, where
$\{ \psi_i (\z) \}$ is the collection of fields
appearing in operator expansion~\eqref{e2.13}.

\section{Generalized vertex algebras}\label{s3}

\subsection{Parafermion fields}\label{s3.1}

Now we will introduce a grading on the vector space $V$ and on the 
space of fields, which will allow us to
make the notion of locality more concrete.

Let us fix a (not necessarily finite or discrete) abelian group $Q$,
and let us assume that our vector space $V$ is \emph{$Q$-graded}:
$V = \bigoplus_{\al\in Q} V_\al$. 
In addition, assume we are given a 
symmetric bilinear map
$\De\colon Q \times Q \to \CC/\ZZ$.
As before, we will identify cosets $\Ga \in \CC/\ZZ$ 
with subsets of $\CC$ of the form $d+\ZZ$ for some $d\in\CC$.
Note that, although $\De(\al,\be)$ is defined mod $\ZZ$, 
$e^{ -2\pi\i \De(\al,\be) }$
is a well-defined complex number.

\skp
\Definition{ \ref{s3.1}}
{\rm (a)}
A \emph{parafermion field} (or just a \emph{field} for short)
of \emph{degree} $\al\in Q$ is a formal series 
$a(\z) \in (\End V) \llb \z,\z^\CC \rrb$
with the property that 
\begin{equation}\label{e3.2}
a(\z) b \in V_{\al+\be} \llb \z \rrb \z^{-\De(\al,\be)}
\qquad \text{for} \quad b \in V_\be \,.
\end{equation}
We denote by $\QF_\al(V)$ the vector space of all fields of degree $\al$,
and by $\QF(V) := \bigoplus_{\al\in Q} \QF_\al(V)$
the vector space of all \emph{fields} on~$V$.

\medskip

{\rm (b)}
Two parafermion fields $a(\z)$, $b(\z)$ of degrees $\al$ and $\be$,
respectively, are called \emph{mutually local\/} if
(see Eqs.\ \eqref{e2.3}, \eqref{e2.4})
\begin{equation}\label{e3.9}
\io_{\z,\w} (\z-\w)^{N} \, a(\z) b(\w)
= \ett(\al,\be) \, \io_{\w,\z} (\z-\w)^{N} \, b(\w) a(\z)
\end{equation}
for some $N \in \De(\al,\be)$ and $\ett(\al,\be) \in \CCC$.
\skp

It is easy to see that Eq.\ \eqref{e3.9} forces $N \in \De(\al,\be)$
and the following relation:
\begin{equation}\label{e3.1}
\ett(\al,\be) \, \ett(\be,\al) = e^{ -2\pi\i \De(\al,\be) } 
\,, \qquad \al,\be \in Q \,.
\end{equation}
Then the above definition of locality is symmetric with respect to
switching $a(\z)$ and $b(\z)$, due to Eqs. \eqref{e2.3}, \eqref{e2.4},
\eqref{e3.1}. The notion of locality can be extended to not necessarily
homogeneous fields (i.e., of fixed degree) by requiring that all their
homogeneous components be mutually local.

\skp
\Example{ \ref{s3.1}}
For $Q = \ZZ/2\ZZ = \{\bar0,\bar1\}$, we can think of 
$V = V_{\bar0} \oplus V_{\bar1}$ as a \emph{superspace}.
Then the choice $\De(\al,\be) = \ZZ$, $\ett(\al,\be) = (-1)^{\al\be}$
corresponds to the usual locality for fields on $V$ 
(see \cite{G,DL,K,Li1}).

\skp
\Remark{ \ref{s3.1}}
Our notation coincides
with that 
in \cite[Chapter 0]{FFR}, except that in \cite{FFR} the abelian
group $Q$ is denoted by $\Ga$ and is assumed finite.
Let us compare our notation to that in \cite{DL}. 
What is denoted by
$(\al,\be)$ in \cite{DL} corresponds to our $\De(\al,\be)$;
however, in \cite{DL} all $(\al,\be)$ belong to a finite subgroup
of $\CC / 2\ZZ$, while $\De(\al,\be)$ take values in $\CC / \ZZ$.
The function $c(\al,\be)$ from \cite{DL}
coincides with $\ett(\al,\be) \, e^{ \pi\i \De(\al,\be) }$ in our
notation; then the condition 
$c(\al,\be) = c(\be,\al)^{-1}$ 
is equivalent to Eq.\ \eqref{e3.1}
above. We prefer to work with $\ett(\al,\be)$ instead of $c(\al,\be)$
because $e^{ \pi\i \De(\al,\be) }$ is not well defined
when $\De(\al,\be)$ is defined mod $\ZZ$.

\skp

We are going to write fields in the conventional form
\begin{equation}\label{e3.3}
a(\z) = \sum_{n\in\CC} a_{(n)} \z^{-n-1} \,, 
\qquad a_{(n)} \in \End V \,.
\end{equation}
Then Eq.\ \eqref{e3.2} is equivalent to the following 
conditions:
\begin{equation*}
a_{(n)} b \in V_{\al+\be} \quad\text{for}\quad
a(\z) \in \QF_\al(V) \,, \;\; b\in V_\be 
\end{equation*}
and
\begin{equation}\label{e3.5}
a_{(n)} b = 0 \quad\text{if}\quad 
n \notin \De(\al,\be)
 \quad\text{or}\quad \Re\, n \gg 0 \,.
\end{equation}
The coefficients $a_{(n)}$ in Eq.\ \eqref{e3.3} are called
\emph{modes} of $a(\z)$.

Assume, in addition, that we are given the 
\emph{vacuum vector} $\vac \in V_0$
and \emph{translation operator\/} $T \in \End V$
preserving the $Q$-grading of $V$ and such that $T\vac=0$.
As before, a field $a(\z)$ is called
\emph{translation covariant\/} if
\begin{equation*}
T a(\z) - a(\z) \, T = \di_{\z} a(\z) \,.
\end{equation*}
In this case, \prref{s2.2}(c) implies that
$a(\z) \vac \in V \llb\z\rrb$. 
On the other hand, the bilinearity of $\De$ implies
$\De(\al,0) = \ZZ$ for all $\al\in Q$.
Thus, every translation covariant parafermion field 
$a(\z)$ of degree $\al$ satisfies
$a(\z) \vac \in V_\al \llb\z\rrb$. 

\subsection{Completeness and state-field correspondence}\label{s3.2}

As in the previous subsection, 
let $V$ be a $Q$-graded vector space endowed with a translation operator
$T$ and a vacuum vector $\vac$.

\skp
\Definition{ \ref{s3.2}}
A system of fields $\{ \phi_{i} (\z) \}_{ i \in \I }$
is called \emph{local} iff
$\phi_{i} (\z)$ and $\phi_{j} (\z)$ are mutually local
for every $i,j \in \I$.
The system $\{\phi_{i} (\z)\}$ is called
\emph{translation covariant}
iff every $\phi_{i} (\z)$ is translation covariant.
Finally, the system $\{\phi_{i} (\z) \}$
is called \emph{complete} 
iff the coefficients of all formal series
$\phi_{i_1} (\z_1) \dotsm \phi_{i_n} (\z_n) \vac$
($n\in\NN$) together with $\vac$ span the whole vector space $V$.
\skp

The next result shows that, given a
translation covariant, local and complete system of fields on $V$,
one can extend it uniquely to a \emph{state-field correspondence}.

\skp
\Theorem{ \ref{s3.2}}
{\it{
Let\ $\{ \phi_{i} (\z) \}_{ i \in \I }$ be a
translation covariant, local and complete sys\-tem of fields on $V$.
Then for every $a \in V$ there exists a unique field, denoted as
$Y (a,\z)$, which is translation covariant, local with respect to
all\/ $\phi_{i} (\z)$, and such that\/
$Y (a,\z) \vac |_{\z = 0} = a$.
}}

\skp

The \emph{proof\/} follows closely that of Theorem 4.2 from \cite{BN}.
Let us consider the vector space $F$ of all translation covariant
fields that are local with respect to $\phi_{i} (\z)$
for all $i \in \I$.
By \prref{s2.2}(c), there is a well-defined linear map
\begin{equation*}
\Phi\colon F \to  V \,, \qquad
\chi (\z) \mapsto \chi (\z) \vac \big|_{\z = 0} \,.
\end{equation*}
Our goal is to show that this map is an isomorphism of vector spaces.
Without loss of generality, we will assume that the system of fields
$\{ \phi_{i} (\z) \}$ is \emph{homogeneous},
i.e., that every field $\phi_{i} (\z)$ has a certain degree~$\al_i$.

Consider for every fixed
$m \in\NN$ and $i_1,\dots,i_m \in \I$
the $m$-field
\begin{equation*}\label{e3.11}
A (\z_1,\dots,\z_m) :=
\Bigl(
\prod_{1 \le k < l \le m} \io_{\z_k,\z_l} (\z_k-\z_l)^{N_{kl}}
\Bigr) \,
\phi_{i_1} (\z_1) \dotsm \phi_{i_m} (\z_m) \,,
\end{equation*}
where $N_{kl}$ are the numbers fulfilling locality
condition \eqref{e3.9} for $\phi_{i_k}(\z)$ and $\phi_{i_l}(\z)$.
Note that $A$ is a translation covariant $m$-field
and is local with respect to all $\phi_{i} (\z)$,
due to \prref{s2.2}(e).
Then all fields $\psi_i(\z)$ appearing in the operator expansion of $A$
are contained in $F$ (see Eq.\ \eqref{e2.13}).
By \reref{s2.3}, all coefficients of
$A (\z_1,\dots,\z_m) \vac$ belong to
the image of the above map $\Phi$.

On the other hand, the product
$\phi_{i_1} (\z_1) \dotsm \phi_{i_m} (\z_m) \vac$
belongs to the space 
$V \llb\z_1\rrb \z_1^\CC \dotsm \llb\z_m\rrb \z_m^\CC$,
which is a module over the algebra
$\CC \llb\z_1\rrb \z_1^\CC \dotsm \llb\z_m\rrb \z_m^\CC$
(see \seref{s2.1} and \ref{s2.2}).
For $k<l$ each $\io_{\z_k,\z_l} (\z_k-\z_l)^{N_{kl}}$
is invertible in the latter algebra, and we have
\begin{equation*}\label{e3.12}
\phi_{i_1} (\z_1) \dotsm \phi_{i_m} (\z_m) \vac
= \Bigl(
\prod_{1 \le k < l \le m} \io_{\z_k,\z_l} (\z_k-\z_l)^{-N_{kl}}
\Bigr) \,
A (\z_1,\dots,\z_m) \vac \,.
\end{equation*}
This implies that every coefficient of 
$\phi_{i_1} (\z_1) \dotsm \phi_{i_m} (\z_m) \vac$
can be expressed as a linear combination of coefficients of
$A (\z_1,\dots,\z_m) \vac$, and hence belongs to
the image of $\Phi$.
Then completeness of the system $\{ \phi_{i} (\z) \}$ 
implies that the map $\Phi$ is surjective.

To prove injectivity, first notice that translation covariance implies
$\chi (\z) \vac = e^{\z \spr T} 
( \chi (\w) \vac |_{\w = 0} )$.
In particular, every element $\chi (\z)$ of the kernel of
$\Phi$ satisfies $\chi (\z) \vac = 0$.
Without loss of generality we can assume that 
$\chi (\z)$ is a field of degree $\al$ for some $\al\in Q$.
Then, by \prref{s2.2}(e), $\chi (\z)$ is local with respect to the 
field $A (\z_1,\dots,\z_m)$ (see Eq.\ \eqref{e2.12}):
\begin{equation*}
\begin{split}
\Bigl( \prod_{l=1}^m \,
\io_{\z,\z_l} & (\z-\z_l)^{N_l} \Bigr) \,
\chi (\z) A (\z_1,\dots,\z_m)
\\
&= \ett \,
\Bigl( \prod_{l=1}^m \,
\io_{\z_l,\z} (\z-\z_l)^{N_l} \Bigr)
A (\z_1,\dots,\z_m) \chi (\z)
\end{split}
\end{equation*}
for some complex numbers $N_l$ and $\ett$.
If we apply both sides of this equation to the vacuum $\vac$,
the right-hand side becomes zero. Then the same argument as above
will give that $\chi (\z)$ vanishes on all coefficients of
$A (\z_1,\dots,\z_m) \vac$. Since these coefficients span $V$,
we obtain that $\chi (\z)=0$.
This completes the proof of the theorem.\qed

\skp

{}From the proof of the above theorem
one can deduce the following generalization of
the Uniqueness Theorem from \cite{G,K}.

\skp
\Corollary{ \ref{s3.2}}
{\it{
Let\ $\{ \phi_{i} (\z) \}$ be a local and complete
system of fields on $V$, and let\/
$B$ be an $n$-field that is
local with respect to all\/ $\phi_{i} (\z)$. Then
if\/ $B(\z_1,\dots,\z_n)$ annihilates the vacuum vector, 
it is identically zero on\/ $V$.
}}

\skp
\Remark{ \ref{s3.2}}
(cf.\ \cite{DK}).
It follows from the above Theorem that the collection 
of fields $\{ Y(a,\z) \}_{a \in V}$ 
coincides with the set of all parafermion fields
that are translation covariant and local with respect to 
all~$\phi_{i} (\z)$. 

\subsection{Definition of generalized vertex algebra}\label{s3.3}

Motivated by \thref{s3.2}, we introduce the notion 
of a generalized vertex algebra in terms of locality
in the spirit of \cite{K}
(cf.\ \cite{G,DL,Li1,GL}).
We will show in \seref{s3.5} below that our definition
is a generalization of the original ones from \cite{FFR,DL,M}.
Fix an abelian group $Q$ and a symmetric bilinear map
$\De\colon Q \times Q \to \CC/\ZZ$.

\skp
\Definition{ \ref{s3.3}}
A \emph{generalized vertex algebra} consists of the following 
\textbf{data}:

\medskip

(\emph{space of states}) \,
a $Q$-graded vector space $V = \bigoplus_{\al\in Q} V_\al\,$;

\medskip

(\emph{vacuum vector}) \,
a vector $\vac \in V_0\,$;

\medskip

(\emph{translation operator}) \,
a $Q$-grading preserving endomorphism $T \in\End V$;

\medskip

(\emph{state-field correspon\-dence}) \,
a $Q$-grading preserving 
linear map\/
$Y \colon V \to \QF(V)$, \, 
$a \mapsto Y(a,\z) = \sum_{n\in\CC} \, a_{(n)} \z^{-n-1}$ \,
from $V$ to the space of parafermion fields on $V$,

\medskip\noindent
subject to the following \textbf{axioms}
($a \in V_\al$, $b \in V_\be$):

\medskip

(\emph{vacuum axiom}) \,
$T\vac=0$, \,
$Y (a,\z) \vac |_{\z = 0} = a\,$;

\medskip

(\emph{translation covariance}) \,
$[T, Y (a,\z)] = \di_{\z} Y (a,\z)\,$;

\medskip

(\emph{locality}) \,
$\io_{\z,\w} (\z-\w)^{N} \, Y(a,\z) Y(b,\w)$
\\ \hspace*{90pt}
$= \ett(\al,\be) \, \io_{\w,\z} (\z-\w)^{N} \, Y(b,\w) Y(a,\z)$

\smallskip\noindent
for some $N \in \De(\al,\be)$ and a bimultiplicative function 
$\ett\colon Q \times Q \to \CCC$
such that 
$\ett(\al,\be) \, \ett(\be,\al) = e^{ -2\pi\i \De(\al,\be) }$.

\medskip

A \emph{homomorphism} of generalized vertex algebras is a linear map 
$f \colon V \to W$, preserving the $Q$-grading, mapping the vacuum vector
to the vacuum vector, intertwining the actions of the translation operators,
and preserving all $n$-th products:
$f(a_{(n)} b) = f(a)_{(n)} f(b)$. 

\skp
\Remark{ \ref{s3.3}}
(cf.\ \cite{DL}).
One can define a more general notion of a homomorphism 
as a pair $(f,\ph)$, where $f \colon V \to W$ is a linear map
and $\ph\colon Q \to Q$ is a group homomorphism such that
$\ett(\ph\al,\ph\be) = \ett(\al,\be)$. The map $f$ should have the
same properties as in the above Definition, except that instead of preserving
the $Q$-grading it satisfies $f(V_\al) \subseteq W_{\ph\al}$.

\skp

As a corollary of \thref{s3.2}, we obtain a generalization of
the Existence Theorem from \cite{FKRW,K} 
(see also \cite{Li1,GL,LL}).

\skp
\Corollary{ \ref{s3.3}}
{\it{
Every translation covariant, local and complete
system of para\-fermion fields\/ $\{\phi_{i} (\z) \}$ on
a\/ $Q$-graded vector space $V$
generates on $V$ a unique structure 
of a generalized vertex algebra.
}}

\skp
\Examples{ \ref{s3.3}}
(a)
For $Q=\{0\}$ 
the notion of a generalized vertex algebra coincides with that
of a \emph{vertex algebra}.

(b)
For $Q = \ZZ/2\ZZ$ and
$\ett(\al,\be) = (-1)^{\al\be}$, 
a generalized vertex algebra is the same as
a \emph{vertex superalgebra} (cf.\ \exref{s3.1}).

(c)
Let $V$ be a generalized vertex algebra such that
$\ett(\al,\be) = (-1)^{p(\al) p(\be)}$ for
some homomorphism $p \colon Q \to \ZZ / 2\ZZ$.
Define a structure of a vector superspace on $V$ by letting the parity
of $a \in V_\al$ be $p(\al)$. Then $V$ is a 
\emph{vertex superalgebra} and it is
\emph{$Q$-graded}, i.e., 
$a_{(n)} b \in V_{\al+\be}$ for all $a \in V_\al$, $b \in V_\be$, 
$n\in\ZZ$.

\skp

We will give less obvious examples of generalized vertex algebras
in \seref{s5} below; see \cite{FFR,DL,GL} for additional examples.
{}From now on, we will often use the notation
$a(\z) \equiv Y(a,\z)$,
and we will denote the modes  of the field
$Y(a,\z)$ by $a_{(n)}$ as in Eq.\ \eqref{e3.3}.
As already noticed before (see the proof of \thref{s3.2}),
translation covariance implies that
\begin{equation}\label{e3.14}
Y(a,\z) \vac = e^{\z \spr T} a \,, \qquad a \in V \,,
\end{equation}
and in particular $Ta = a_{(-2)} \vac$.
This shows that the translation operator\/ $T$ is uniquely determined
by the state-field correspon\-dence $Y$.
Another property, which follows from \thref{s3.2}, is that
$Y(Ta,\z) = \di_\z Y(a,\z)$.
We also have the \emph{skew-symmetry} relation \cite{FFR,GL}:
\begin{equation}\label{e3.15s}
Y(a,\z) b = 
\ett(\al,\be) \, 
e^{\z \spr T} \bigl( Y(b,e^{\pi\i} \z) a \bigr) \,,
\qquad a \in V_\al \,,  \; b \in V_\be \,.
\end{equation}
This can be proved by applying both sides of Eq.\ \eqref{e3.9} to $\vac$,
using \eqref{e3.14} and translation covariance, 
and setting $\w=0$.

\subsection{Action of the cohomology group $\coh^2(Q,\CCC)$}\label{s3.6}

In this subsection we show that the second cohomology group $\coh^2(Q,\CCC)$
acts naturally on the isomorphism classes of generalized vertex algebras
with fixed $Q$ and $\De$. To the best of our knowledge, this action is new;
however, the idea of modifying the vertex operators by a $2$-cocycle
goes back to \cite{FK}.

Let $Q$ be a fixed abelian group.
It is well known that $\coh^2(Q,\CCC)$ parameterizes the
\emph{central extensions} of the group $Q$, i.e., 
exact sequences of group homomorphisms
\begin{equation}\label{e4.1}
1 \to \CCC \to \wti Q \to Q \to 1
\end{equation}
such that the image of $\CCC$ is central in $\wti Q$,
up to equivalence.
We start this subsection by reviewing the theory of such
central extensions (see e.g.\ \cite[Chapter IV]{Br}).
This material will be used here
and in \seref{s5.2} below. As before, we will write
the group operation in $Q$ additively, while the one in $\wti Q$
will be written multiplicatively.

Given a section $e \colon Q \to \wti Q$ 
of extension \eqref{e4.1}, one can identify
$\wti Q$ as a set with $\CCC \times Q$, so that
$e^\al := e(\al) \in \wti Q$ is identified with 
$(1,\al)$ for $\al \in Q$, and the embedding
$\CCC \to \wti Q$ is given by $c \mapsto (c,0)$.
Then all elements of $\wti Q$ have the form
$c \, e^\al \equiv (c,\al)$.
The product in $\wti Q$ gives rise to a map
$\ep\colon Q \times Q \to \CCC$ such that
\begin{equation}\label{e4.5e}
e^\al e^\be = \ep(\al,\be) \, e^{\al+\be} \,,
\qquad \al,\be\in Q \,.
\end{equation}
Associativity and unit properties of this product are equivalent to the
condition that $\ep$ is a (normalized)
\emph{$2$-cocycle}, i.e.,
\begin{equation}\label{e4.6}
\begin{split}
\ep(\al,\be) \, \ep(\al+\be, \ga) &=
\ep(\be, \ga) \, \ep(\al, \be+\ga) \,,
\\
\ep(\al,0) &= \ep(0,\al) = 1
\,.
\end{split}
\end{equation}
If we choose a different section $e' \colon Q \to \wti Q$,
then it is related to the section $e$ by
$e'(\al) = \rho(\al) \, e^\al$
for some function $\rho \colon Q \to \CCC$.
The $2$-cocycle corresponding to the section $e'$ is
\begin{equation}\label{e4.7}
\ep'(\al,\be) = 
\ep(\al,\be) \, \rho(\al) \, \rho(\be) \, \rho(\al+\be)^{-1}
\,,
\end{equation}
so it belongs to the same \emph{cohomology class} as the $2$-cocycle $\ep$.
The group (with respect to multiplication of functions) of all such
cohomology classes is the \emph{cohomology group} $\coh^2(Q,\CCC)$.
This gives a one-to-one correspondence between
equivalence classes of
central extensions \eqref{e4.1} and elements of $\coh^2(Q,\CCC)$.

For a given section $e \colon Q \to \wti Q$ as above, one finds that
\begin{equation*}
e^\al e^\be = \oom(\al,\be) \, e^\be e^\al \,,
\qquad \al,\be\in Q \,,
\end{equation*}
where
\begin{equation*}
\oom(\al,\be) = \ep(\al,\be) \, \ep(\be,\al)^{-1} \,.
\end{equation*}
Clearly, $\oom$ does not depend on the choice of section,
so it depends only on the cohomology class of $\ep$.
Thus $\oom$ is an invariant of the central extension \eqref{e4.1},
which we will call its \emph{canonical invariant}.
Introduce the group (with respect to multiplication)
$\OM(Q)$ of all bimultiplicative maps
$\oom\colon Q \times Q \to \CCC$
such that $\oom(\al,\al) = 1$ for all $\al\in Q$.
Then the canonical invariant $\oom$ of any
central extension \eqref{e4.1} belongs to $\OM(Q)$,
and we get a group homomorphim from $\coh^2(Q,\CCC)$ to $\OM(Q)$.
The next result is well known.\footnote{We thank Ofer Gabber for a discussion
on this.}

\skp
\Lemma{ \ref{s3.6}}
{\it{
The above homomorphim from\/ $\coh^2(Q,\CCC)$ to\/ $\OM(Q)$ is 
an isomorphism of groups. In particular, for every\/
$\oom\in\OM(Q)$ there exists a central extension \eqref{e4.1} 
with a canonical invariant\/ $\oom$.
}}

\skp
\Proof
It follows from, e.g., \cite[\S V.6, Exercise 5]{Br} that
the homomorphim $\coh^2(Q,\CCC)\to\OM(Q)$
is surjective and the kernel is isomorphic to
$\Ext_\ZZ(Q,\CCC)$. But $\CCC$ is isomorphic to $\CC/\ZZ$
as a $\ZZ$-module. The latter is divisible, and hence is an
injective $\ZZ$-module (see \cite[\S III.4]{Br}). Therefore,  
$\Ext_\ZZ(Q,\CCC) = \{0\}$.\qed


\skp

Now let $V$ be a generalized vertex algebra as in \deref{s3.3}.
Given a $2$-cocycle $\ep\colon Q \times Q \to \CCC$, 
we modify the state-field correspon\-dence $Y$ by defining
\begin{equation}\label{e4.9}
Y^\ep (a,\z) b := \ep(\al,\be) \, Y(a,\z) b \,,
\qquad a \in V_\al \,, \; b \in V_\be \,.
\end{equation}
Then $Y^\ep$ satisfies the vacuum and translation covariance axioms
(with the same $\vac$ and $T$) and the locality axiom with
$\ett(\al,\be)$ replaced by
\begin{equation*}
\ett^\ep(\al,\be) := 
\ett(\al,\be) \, \ep(\be,\ga) \, \ep(\al,\be+\ga) 
\, \ep(\al,\ga)^{-1} \, \ep(\be,\al+\ga)^{-1}
\,.
\end{equation*}
Using Eq.\ \eqref{e4.6}, this can be simplified to
\begin{equation}\label{e4.11}
\ett^\ep(\al,\be) := \ett(\al,\be) \, 
\ep(\al,\be) \, \ep(\be,\al)^{-1} 
= \ett(\al,\be) \, \oom(\al,\be) \,,
\end{equation}
where $\oom$ is the canonical invariant of the central extension
defined by $\ep$.
Hence $\ett^\ep$ is bimultiplicative and satisfies Eq.\ \eqref{e3.1} 
with the same $\De$.
Therefore, we obtain a new structure of a generalized vertex algebra
on $V$ with the same $Q$-grading, $\De$, $\vac$ and $T$ but with
modified $Y$ and $\ett$.

\skp
\Definition{ \ref{s3.6}}
The above-defined 
\emph{modified} 
generalized vertex algebra will be denoted as $V^\ep$.
We say that two generalized vertex algebras $V$ and $W$
(for the same $Q$, $\De$) are \emph{equivalent}
if $W$ is isomorphic to $V^\ep$ for some $2$-cocycle $\ep$.

\skp

This is an equivalence relation because $V^{\ep_1 \ep_2}$ 
is naturally isomorphic to $(V^{\ep_1})^{\ep_2}$.
Note also that, for every homomorphism
$f\colon V\to W$ of generalized vertex algebras,
the same map $f$ is a homomorphism
from $V^\ep$ to $W^\ep$ (since $f$ preserves the $Q$-grading).
More generally, let $(f,\ph) \colon V \to W$ be a homomorphism
in the sense of \reref{s3.3}. Then for every $2$-cocycle $\ep$,
$\ep'(\al,\be) := \ep(\ph\al,\ph\be)$ is also a $2$-cocycle,
and $(f,\ph) \colon V^{\ep'} \to W^{\ep}$ is a homomorphism.

It is obvious from Eq.\ \eqref{e4.11} that 
$\ett^\ep$ depends only on the cohomology
class $[\ep]$ of $\ep$ in $\coh^2(Q,\CCC)$.
Given $\ep' \in [\ep]$
(see Eq.\ \eqref{e4.7}), we define an endomorphism $\rho$ of $V$
by $\rho(a) := \rho(\al) a$ for $a \in V_\al$. Then \eqref{e4.7} 
implies that
$\rho \colon V^{\ep'} \to V^{\ep}$ 
is an isomorphism of generalized vertex algebras.
We thus obtain an action of the cohomology
group $\coh^2(Q,\CCC)$ on the isomorphism classes of 
generalized vertex algebras.

We finish this subsection with an important application of the above results.

\skp
\Proposition{ \ref{s3.6}}
{\it{
Let\/ $V$ be a generalized vertex algebra for which 
the map\/ $\De\colon Q \times Q \to \CC / \ZZ$
is trivial, so that\/
$\ett(\al,\al)^2 = 1$ for all\/ $\al\in Q$.
Then $V$ is equivalent to a $Q$-graded vertex superalgebra,
in which the parity of\/ $a \in V_\al$ is\/
$p(\al) \in \ZZ / 2\ZZ$
where\/ $(-1)^{p(\al)} = \ett(\al,\al)$.
}}

\skp
\Proof
The function $\oom(\al,\be) := (-1)^{p(\al) p(\be)} \ett(\al,\be)^{-1}$
is bimultiplicative and $\oom(\al,\al) = 1$ for all $\al$.
By \leref{s3.6}, there exists a central extension \eqref{e4.1}
whose canonical invariant is $\oom$. Let
$\ep$ be a $2$-cocycle defining this central extension.
Then $\ett^\ep(\al,\be) = (-1)^{p(\al) p(\be)}$,
and hence $V^\ep$ is a $Q$-graded vertex superalgebra
(see Example~\ref{s3.3}(c)).\qed

\subsection{Operator product expansion and associativity}\label{s3.4}

In this subsection we investigate the operator product expansion
of two parafermion fields in a generalized vertex algebra $V$
(cf.\ \seref{s2.3}).

Let $a \in V_\al$, $b \in V_\be$ for some $\al,\be \in Q$,
and as before let $a(\z) \equiv Y(a,\z)$ for short.
Denote by $F_{a,b} (\z,\w)$ the left-hand side of Eq.\ \eqref{e3.9}.
By \prref{s2.2}(e), $F_{a,b} (\z,\w)$
is a translation covariant bilocal field, which is
local with respect to $c(\z)$ for all $c \in V$.
Then the \emph{operator product expansion} (OPE) 
of $a(\z)$ and $b(\w)$ is defined by Eq.\ \eqref{e2.13} for 
$A(\z_1,\z_2) = F_{a,b} (\z_1,\z_2)$, i.e.,
it is the expansion
\begin{equation}\label{e3.30}
\iota_{\z,\w} \,
F_{a,b} (\z+\w,\z)
= e^{ \w \spr \Di_{\z_1} } F_{a,b} (\z_1,\z_2)
\big|_{ \z_1 = \z_2 = \z  } \,.
\end{equation}
All coefficients in front of powers of $\w$ in this expansion
are translation covariant fields and are local with respect to
$c(\z)$ for all $c \in V$ (see \prref{s2.2}).
Therefore, by \thref{s3.2},
they must be themselves of the form $c(\z)$ for some $c \in V$.

\skp
\Lemma{ \ref{s3.4}}
{\it{
For every two elements\/ $a,b$ 
in a generalized vertex algebra, we have{\rm:}
\begin{equation}\label{e3.29}
\frac1{k!} \,
\di_{\z_1}^k F_{a,b} (\z_1,\z_2) \big|_{\z_1=\z_2=\z}
= Y(a_{(N-1-k)} b,\z)
\,, \qquad k \in \NNN \,,
\end{equation}
where\/ $N$ is from Eq.\ \eqref{e3.9} and\/ $F_{a,b}$ 
denotes the left-hand side of\/ Eq.~\eqref{e3.9}.
}}

\skp
\Proof
By the above discussion, it is enough to establish Eq.\ \eqref{e3.29}
after both sides are applied to $\vac$ and $\z$ is set to $0$.
But, by \prref{s2.2}(c), 
$F_{a,b} (\z_1,\z_2) \vac \in V\llb\z_1,\z_2\rrb$.
Then
$F_{a,b} (\z_1,\z_2) \vac |_{\z_2=0} = \z_1^N a(\z_1) b \in V\llb\z_1\rrb$;
hence applying to it $\frac1{k!} \, \di_{\z_1}^k$ 
and evaluating it at $\z_1=0$ gives the coefficient in front of 
$\z_1^k$ in $\z_1^N a(\z_1) b$,
which is exactly $a_{(N-1-k)} b$.\qed

\skp

Let us point out that as a corollary of the above proof, we have
$\z^N a(\z) b \in V\llb\z\rrb$. 
Combining Eqs.\ \eqref{e3.30} and \eqref{e3.29}, we obtain that
the OPE of $a(\z)$ and $b(\w)$ is given by:
\begin{equation}\label{e3.31}
\iota_{\z,\w} \,
F_{a,b} (\z+\w,\z) = \w^N \, Y(a(\w) b, \z) \,.
\end{equation}
Denote by $N(a,b)$ the minimal $N \in \De(\al,\be)$
(i.e., with minimal real part) fulfilling
locality condition \eqref{e3.9}.
Now we can prove the following \emph{formal associativity}
and \emph{commutativity} relations generalizing those of \cite{FB,LL}
(see also \cite{DL,BN}).

\skp
\Theorem{ \ref{s3.4}}
{\it{
Let\/ 
$a \in V_\al$, $b \in V_\be$, $c \in V_\ga$ $(\al,\be,\ga \in Q)$
be three elements in a generalized vertex algebra\/ $V$.
Then there exists a localized formal series
\begin{equation*}
\Y_{a,b,c} (\z,\w) = 
\psi_{a,b,c} (\z,\w) \,
\z^{-N (a,c)} \, \w^{-N (b,c)} \, (\z-\w)^{-N (a,b)}
\,,
\end{equation*}
where\/
$\psi_{a,b,c} (\z,\w) \in V \llb \z,\w \rrb$,
with the properties that
}}
\begin{align*}
a (\z)  b (\w)  c  &=
\iota_{\z,\w} \, \Y_{a,b,c} (\z,\w) \,,
\raisebox{10pt}{} \\ 
\ett(\al,\be) \,
b (\w)  a (\z)  c  &=
\iota_{\w,\z} \, \Y_{a,b,c} (\z,\w) ,
\raisebox{10pt}{} \\ 
Y( a (\w)  b, \z)  c  &=
\iota_{\z,\w} \, \Y_{a,b,c} (\z+\w,\z) .
\end{align*}
\Proof
We set $N=N(a,b)$ in Eq.\ \eqref{e3.9}, and define
\begin{equation*}
\Y_{a,b,c} (\z,\w) := (\z-\w)^{-N (a,b)} \, F_{a,b} (\z,\w) c \,.
\end{equation*}
Then everything follows from Eqs.\ \eqref{e3.9}, \eqref{e3.31}
and the above observation that 
$\z^{N(a,b)} \, a(\z) b \in V\llb\z\rrb$.\qed

\skp

{}From this theorem we deduce the following \emph{associativity}
relations (see \cite{DL,M,BK}). 

\skp
\Corollary{ \ref{s3.4}}
{\it{
For every three elements\/ 
$a \in V_\al$, $b \in V_\be$, $c \in V_\ga$ $(\al,\be,\ga \in Q)$
in a generalized vertex algebra, and for
each\/ $L \in \De(\al,\ga)$,
$L \ge N(a,c)$, we have{\rm:} 
\begin{equation}\label{e3.38}
\iota_{\z,\w} (\z+\w)^L \; \iota_{\z,\w} a(\z+\w) b(\w) c
=
\iota_{\w,\z} (\z+\w)^L \, Y( a(\z) b, \w) c \,,
\end{equation}
\begin{equation*}
\z^L \, a(\z) b(\w) c
=
\Bigl[ \io_{\u,\z-\w} (\u+\z-\w)^L \, \io_{\z,\w}
Y( a (\z-\w) b, \u) c
\Bigr]_{\u = \w} \,.
\end{equation*}
The expression under the substitution\/ $\u=\w$ belongs to the space\/
$\u^{-N(b,c)}$ $\times \,
\io_{\z,\w} (\z-\w)^{-N(a,b)} \, V \llb\z,\w,\u\rrb$,
and hence the substitution makes sense.
}}

\skp
The \emph{proof\/} is the same as that of Proposition 5.4 from \cite{BN},
utilizing Taylor's formula \eqref{e2.15}.\qed

\subsection{Jacobi identity and Borcherds identity}\label{s3.5}

In this subsection we will show that every generalized vertex algebra $V$
obeys the Jacobi identity of \cite{FFR,DL,M}, and hence our definition 
is a generalization of the ones from \cite{FFR,DL,M}.

First, recall that the \emph{formal delta-function} 
is defined as the formal series 
(see Eqs.\ \eqref{e2.3}, \eqref{e2.4}):
\begin{equation}\label{e3.17}
\de(\z,\w) := 
(\io_{\z,\w}- \io_{\w,\z}) (\z-\w)^{-1}
= \sum_{j\in\ZZ} \z^{-j-1} \w^j \,.
\end{equation}
More generally, for a subset $\Ga\subseteq\CC$, we define 
(see \seref{s2.1})
\begin{equation}\label{e3.18}
\de_{\Ga} (\z,\w) := \sum_{j\in\Ga} \z^{-j-1} \w^j 
\in \CC\llb \z,\z^\Ga, \w,\w^\Ga \rrb \,.
\end{equation}
Note that if $\Ga = d + \ZZ$ for some $d\in\CC$, then
$\de_{\Ga} (\z,\w) = (\w/\z)^{d} \de(\z,\w)$,
and this is independent of the choice of $d$ mod $\ZZ$.

\skp
\Theorem{ \ref{s3.5}}
{\it{
For every three elements\/ 
$a \in V_\al$, $b \in V_\be$, $c \in V_\ga$ $(\al,\be,\ga \in Q)$
in a generalized vertex algebra, and for every\/
$n \in \De(\al,\be)$, we have the following identity{\rm:}
\begin{align}
\notag
Y(a,\z) Y&(b,\w) c \; \io_{\z,\w} (\z-\w)^n
- \ett(\al,\be) \,
Y(b,\w) Y(a,\z) c \; \io_{\w,\z} (\z-\w)^n
\quad
\\ \label{e3.19}
& =
\sum_{j \in\NNN} Y(a_{(n+j)} b,\w) c
\; \di_w^j \de_{\De(\al,\ga)} (\z,\w) / j!
\,. 
\end{align}
The sum in the right-hand side is finite due to Eq.\ \eqref{e3.5}.
}}

\skp
\Proof
As in the proof of Theorem 5.5 from \cite{BN},
we deduce from \thref{s3.4} above that 
\begin{align*}
a(\z) b(\w) c \; & \io_{\z,\w} (\z-\w)^n 
\\
= \z^{-L} \, & \Bigl[ \io_{\u,\z-\w} (\u+\z-\w)^L \, \io_{\z,\w} \, 
(\z-\w)^n  \, Y( a (\z-\w) b, \u) c 
\Bigr]_{\u = \w} \,,
\intertext{and}
\ett(\al,\be) \,
b(\w) & a(\z) c \; \io_{\w,\z} (\z-\w)^n 
\\
= \z^{-L} \, & \Bigl[ \io_{\u,\z-\w} (\u+\z-\w)^L \, \io_{\w,\z} \, 
(\z-\w)^n  \, Y( a (\z-\w) b, \u) c 
\Bigr]_{\u = \w}
\end{align*}
for large enough $L \in \De(\al,\ga)$.
Recall that (see Eq.\ \eqref{e2.3})
\begin{equation*}
\io_{\u,\z-\w} (\u+\z-\w)^L
= \sum_{i\in\NNN} \binom{L}{i} \u^{L-i} \, (\z-\w)^i \,.
\end{equation*}
Take the difference of the above two expressions,
and notice that
\begin{equation*}
\begin{split}
(\io_{\z,\w} & - \io_{\w,\z})
\Bigl( 
\io_{\u,\z-\w} (\u+\z-\w)^L \,
(\z-\w)^n \, a (\z-\w) b 
\Bigr)
\\
& = 
\sum_{i\in\NNN} \sum_{j\in\ZZ} \binom{L}{i} \u^{L-i} \,
(\io_{\z,\w} - \io_{\w,\z}) (\z-\w)^{i-j-1}
\, a_{(n+j)} b \,.
\end{split}
\end{equation*}
Since $(\io_{\z,\w} - \io_{\w,\z}) (\z-\w)^{i-j-1} = 0$
when $i-j-1 \ge 0$, the right-hand side
can be rewritten as follows:
\begin{equation*}
\sum_{j\in\NNN} a_{(n+j)} b \;
\sum_{i=0}^j \binom{L}{i} \u^{L-i} \,
(\io_{\z,\w} - \io_{\w,\z}) (\z-\w)^{i-j-1} \,.
\end{equation*}
The sum over $j\in\NNN$ is in fact finite, because
$a_{(n+j)} b = 0$ for $j \gg 0$ by Eq.~\eqref{e3.5}.
Then both sums over $i$ and $j$ are finite. 

Therefore, the substitution $\u=\w$ makes sense, and we obtain
that the left-hand side of Eq.\ \eqref{e3.19} is equal to
\begin{equation*}
\sum_{j\in\NNN} Y(a_{(n+j)} b, \w) c \;
\sum_{i=0}^j \binom{L}{i} \z^{-L} \, \w^{L-i} \,
(\io_{\z,\w} - \io_{\w,\z}) (\z-\w)^{i-j-1} \,.
\end{equation*}
Now observe that, by Eq.\ \eqref{e3.17},
\begin{equation*}
(\io_{\z,\w}- \io_{\w,\z}) (\z-\w)^{-k-1}
= \di_\w^{(k)} \de(\z,\w)
\,, \qquad k \in\NNN \,,
\end{equation*}
where we use the divided-power notation 
$\di_\w^{(k)} := \di_\z^k / k!$.
This formula and the fact that 
$\di_\w^{(i)} \w^L = \binom{L}{i} \w^{L-i}$,
together with the Leibniz rule, imply that
the sum over $i$ in the above equation is exactly
$\di_\w^{(j)} \bigl( \z^{-L} \, \w^L \, \de(\z,\w) \bigr)$.
Since, by definition, 
$\z^{-L} \w^L \de(\z,\w) = \de_{\De(\al,\ga)} (\z,\w)$,
this completes the proof.\qed

\skp

We are going to rewrite Eq.\ \eqref{e3.19} in terms of modes.
By definition (see \eqref{e3.3}),
the modes of a field $a(\z)$ are obtained from it
by taking residues:
\begin{equation*}
a_{(n)} = \Res_\z \z^n \, a(\z) \,,
\qquad n \in \CC \,,
\end{equation*}
where $\Res_\z$ gives the coefficient in front of $\z^{-1}$.
Now multiply both sides of Eq.\ \eqref{e3.19} by $\z^m \w^k$,
take $\Res_\z \Res_\w$, use expansions \eqref{e2.3}, \eqref{e2.4},
and note that, because of Eq.\ \eqref{e3.5}, only
$m\in\De(\al,\ga)$ and $k\in\De(\be,\ga)$ will give a nonzero result.
We thus obtain the \emph{Borcherds identity}:
\begin{align}
\notag
\sum_{j\in\NNN} (-1)^j & \binom{n}{j} 
\Bigl( 
a_{(m+n-j)}(b_{(k+j)}c)
- \ett(\al,\be) \, e^{\pi\i n} \,
b_{(n+k-j)}(a_{(m+j)}c)
\Bigr)
\\ \label{e3.28}
&=
\sum_{j\in\NNN} \binom{m}{j} (a_{(n+j)}b)_{(m+k-j)}c \,,
\end{align}
for $a \in V_\al$, $b \in V_\be$, $c \in V_\ga$,
$n\in\De(\al,\be)$, $m\in\De(\al,\ga)$, $k\in\De(\be,\ga)$.

\smallskip

It is remarkable that this identity has exactly the same form
for (ordinary) vertex algebras, in which case
$\ett(\al,\be) = 1$ and $n,m,k \in\ZZ$
(see \cite[(4.8.3)]{K}).
Borcherds used special cases of this identity
in his original definition of the notion 
of a vertex algebra in \cite{B}.

\skp
\Remark{ \ref{s3.5}}
Equation \eqref{e3.28} above coincides with \cite[(0.51)]{FFR},
which is equivalent to the ``Jacobi identity'' \cite[(0.46)]{FFR}.
On the other hand,
the collection of identities \eqref{e3.19} for all $n\in\De(\al,\be)$
is equivalent to the ``Jacobi identity'' \cite[(9.14)]{DL}.
Indeed, the expression
$\u^{-1} \de\bigl( \frac{\z-\w}{\u} \bigr)$ 
used in \cite{FLM,DL} coincides in our notation with
$\io_{\z,\w} \de(\z-\w,\u)$, while
$\u^{-1} \de\bigl( \frac{\w-\z}{-\u} \bigr)$ 
from \cite{FLM,DL} is equal to our $\io_{\w,\z} \de(\z-\w,\u)$.
Thus, if we multiply Eq.\ \eqref{e3.19} by $\u^{-n-1}$ and sum
over $n\in\De(\al,\be)$, we obtain \cite[(9.14)]{DL}
with the corresponding changes in notation (cf.\ \reref{s3.1}).

\skp

We will now show that Borcherds identity \eqref{e3.28},
together with a ``partial vacuum axiom'' can be taken as
an equivalent definition of a generalized vertex algebra.

\skp
\Proposition{ \ref{s3.5}}
{\it{
Let\/ $Q$ be a fixed abelian group with a bimultiplicative map\/
$\ett\colon Q \times Q \to \CCC$, as in Sect.\ {\rm\ref{s3.1}}.
Let\/ $V$ be a\/ $Q$-graded vector space endowed with a vector\/ 
$\vac\in V_0$ and with products\/ 
$a_{(n)} b \in V_{\al+\be}$ for\/ $n\in\CC$, $a \in V_\al$, $b \in V_\be$,
satisfying Eqs.\ \eqref{e3.5},
\eqref{e3.28}, and 
\begin{equation*}
\vac_{(n)} a = \de_{n,-1} \, a \,, \qquad
a_{(-1)} \vac = a \,.
\end{equation*}
Then\/ $V$ is a generalized vertex algebra
with\/ $Ta := a_{(-2)} \vac$.
}}

\skp
\Proof
Since Eq.\ \eqref{e3.28} is equivalent to \eqref{e3.19},
it implies locality of $a(\z)$ and $b(\z)$.
The rest of the proof is as in Proposition 4.8(b) from~\cite{K}.\qed

\skp

We refer to \cite{FFR,DL,M,GL} for consequences of the Jacobi identity
(= Borcherds identity),
and we mention only one simple consequence here.
{}From Eq.\ \eqref{e3.19} and the sentence after \leref{s3.4},
we deduce the following well-known fact.

\skp
\Corollary{ \ref{s3.5}}
{\it{
For every two elements\/ 
$a \in V_\al$, $b \in V_\be$, the minimal number $N=N(a,b)$
$($i.e., with minimal real part$)$ fulfilling
locality condition \eqref{e3.9}
is equal to the mini\-mal $N \in \De(\al,\be)$ such that
$\z^N a(\z) b \in V\llb\z\rrb$.
}}

\section{Modules over generalized vertex algebras
and twisted modules over vertex superalgebras}\label{s4}

\subsection{Modules over generalized vertex algebras}\label{s4.1}

We start this subsection by recalling the notion of a module over
a generalized vertex algebra, following \cite{DL}.
Fix an abelian group $Q$, a symmetric bilinear map
$\De\colon Q \times Q \to \CC/\ZZ$, 
and a generalized vertex algebra $V$, as in \seref{s3.3}.
Let $P$ be a \emph{$Q$-set}, i.e., a set with an action of $Q$
on it, which we will write additively. Assume, in addition, that
we are given an extension of $\De$ to a bilinear map 
$\De\colon Q \times P \to \CC/\ZZ$.
For instance, we may take $P=Q$ with the same $\De$, but it will
be convenient to work in the more general framework of $Q$-sets.

For a $P$-graded vector space $M = \bigoplus_{\ga\in P} M_\ga$,
we define the notion of a (parafermion) field on $M$
in the same way as in \deref{s3.1}(a). 
Namely, a \emph{parafermion field} of \emph{degree} $\al\in Q$ 
is a formal series 
$a(\z) \in (\End M) \llb \z,\z^\CC \rrb$
with the property that 
\begin{equation*}
a(\z) c \in M_{\al+\ga} \llb \z \rrb \z^{-\De(\al,\ga)}
\qquad \text{for} \quad c \in M_\ga \,, \; \ga \in P \,.
\end{equation*}
Then the vector space $\QF(M)$ 
of all parafermion fields on $M$ is again $Q$-graded.

\skp
\Definition{ \ref{s4.1}}
A \emph{module} over a generalized vertex algebra $V$ 
(or a $V$-\emph{module}) is 
a $P$-graded vector space $M = \bigoplus_{\ga\in P} M_\ga$
endowed with a $Q$-grading preserving 
linear map
\begin{equation*}
Y \colon V \to \QF(M) \,, \qquad
a \mapsto Y(a,\z) = \sum_{n\in\CC} \, a_{(n)} \z^{-n-1} \,,
\end{equation*}
from $V$ to the space of parafermion fields on $M$,
such that $Y(\vac,\z) = \id$ and
Eq.\ \eqref{e3.19} holds for all 
$a \in V_\al$, $b \in V_\be$, $c \in M_\ga$, $n\in\De(\al,\be)$,
$(\al,\be \in Q$, $\ga \in P)$.

\skp
\Remark{ \ref{s4.1}}
Equivalently, one can define the notion of a $V$-module
in terms of Borcherds identity \eqref{e3.28}.
One can also replace the Borcherds identity (or Jacobi identity)
in the definition of a module by associativity relation \eqref{e3.38}.
This was proved in \cite{Li1} in the case of vertex algebras but 
the same proof works also for generalized vertex algebras, by making use
of skew-symmetry relation \eqref{e3.15s} 
(cf.\ \cite{Li4} and \reref{s4.2} below).

\skp
\Example{ \ref{s4.1}}
Let $V$ be a generalized vertex algebra 
and let $M$ be a $V$-module, as above.
Assume that 
\begin{equation*}
\ett(\al,\be) = (-1)^{p(\al) p(\be)} \,, \quad
\De(\al,\ga) = \ZZ \,, \qquad
\al,\be \in Q , \; \ga \in P \hspace{1pt},
\end{equation*}
for some homomorphism $p \colon Q \to \ZZ / 2\ZZ$.
Then $V$ is a $Q$-graded vertex superalgebra (see \exref{s3.3}(c)),
and $M$ is a \emph{$P$-graded} $V$-module, i.e.,
$a_{(n)} c \in M_{\al+\ga}$
for $a \in V_\al$, $c \in M_\ga$, $n\in\ZZ$.
\skp

The following special situation provides important examples of modules.

\skp
\Proposition{ \ref{s4.1}}
{\it{
Let $P$ be an abelian group endowed with a symmetric bilinear map\/
$\De\colon P \times P \to \CC/\ZZ$, let\/ $V$ be a 
generalized vertex algebra for $P$, and let\/ $Q$ be a subgroup of $P$.
For a coset\/ $\Ga \in P/Q$ $($considered as a subset of $P)$,
introduce the subspace\/
$V_\Ga := \bigoplus_{\ga \in \Ga} V_\ga$ of\/ $V$.
Then\/ $V_Q$ is a generalized vertex algebra for the group\/ $Q$, 
and each\/ $V_\Ga$ is a\/ $V_Q$-module.
}}

\skp
The \emph{proof} is immediate from the definitions.\qed

\skp

{}From this Proposition and the results of \seref{s3.6}, 
we deduce the following corollary, which
can be used to construct modules over
vertex superalgebras (see \seref{s5.2} below).

\skp
\Corollary{ \ref{s4.1}}
{\it{
In the setting of the above Proposition, assume in addition that\/
$P/Q$ is isomorphic to a direct product of cyclic groups and\/
$\De(\al,\ga) = \ZZ$
for all\/ $\al \in Q$, $\ga \in P$. Then there exists
a $2$-cocycle\/ $\ep\colon P \times P \to \CCC$ such that\/
$V_Q^\ep$ is a\/ $Q$-graded vertex superalgebra and each\/
$V_\Ga^\ep$ is a\/ $\Ga$-graded\/ $V_Q^\ep$-module.
}}

\skp
\Proof
By \prref{s3.6} and its proof, there exists an element $\oom \in \Om(Q)$
such that $V_Q^\ep$ is a vertex superalgebra,
for every $2$-cocycle $\ep \colon Q \times Q \to \CCC$ 
with canonical invariant $\oom$.
Let us choose representatives $\ga_i \in P$ for the generators
of the cyclic factors of $P/Q$, where $i$ runs over some
(possibly infinite) index set. Then we extend $\oom$ to a
bimultiplicative map $P \times P \to \CCC$ by letting
$\oom(\ga_i+\al, \ga_j+\be) = \oom(\al,\be)$
for all $i,j$ and all $\al,\be\in Q$.
Then $\oom\in\Om(P)$, and hence it corresponds to a 
$2$-cocycle $\ep\colon P \times P \to \CCC$ (see \leref{s3.6}).
\qed

\subsection{Twisted modules over vertex superalgebras}\label{s4.2}


Let us recall the definition of a twisted module over 
a vertex (super)algebra \cite{FFR,D2}, 
in the formulation given in \cite{DK}.
Fix an additive subgroup $\Ga \subseteq \CC / \ZZ$,
and assume that $V$ is a \emph{$\Ga$-graded} vertex superalgebra,
i.e., $V = \bigoplus_{\al\in \Ga} V_\al$ so that
$a_{(n)} b \in V_{\al+\be}$ for all $a \in V_\al$, $b \in V_\be$, 
$n\in\ZZ$ (cf.\ \exref{s3.3}(c)).

\skp
\Definition{ \ref{s4.2}}
A $\Ga$-\emph{twisted module} over $V$ is a vector space $M$
endowed with a linear map
\begin{equation*}
Y \colon V \to \QF(M) \,, \qquad
a \mapsto Y(a,\z) = \sum_{n\in\CC} \, a_{(n)} \z^{-n-1} \,,
\end{equation*}
from $V$ to the space of parafermion fields on $M$,
such that $Y(\vac,\z) = \id$,
\begin{equation}\label{e4.4}
Y(a,\z) c \in M \llb \z \rrb \z^{-\al} 
\quad\text{for}\quad
a \in V_\al \,, \; c \in M \,, \; \al \in \Ga \,,
\end{equation}
and the following identity holds (cf.\ Eq.\ \eqref{e3.18}):
\begin{align}
\notag
Y(a,\z) Y&(b,\w) c \; \io_{\z,\w} (\z \! - \! \w)^n
- (-1)^{ p(a) p(b) } \,
Y(b,\w) Y(a,\z) c \; \io_{\w,\z} (\z \! - \! \w)^n
\quad
\\ \label{e4.5}
& =
\sum_{j \in\NNN} Y(a_{(n+j)} b,\w) c
\; \di_w^j \de_\al (\z,\w) / j!
\end{align}
for all $a \in V_\al$, $b \in V$, $c \in M$, 
$n\in\ZZ$, $\al \in \Ga$,
where $a$ and $b$ have parities $p(a)$ and $p(b)$, respectively.

\skp
\Example{ \ref{s4.2}}
Let $\si$ be a \emph{diagonalizable} automorphism of $V$,
and let $\Ga$ be a subgroup of $\CC / \ZZ$ such that
all eigenvalues of $\si$ have the form $e^{2\pi\i\al}$ with $\al\in\Ga$.
For example, when $\si^N = 1$, we can take $\Ga = \tfrac1N \ZZ / \ZZ$.
Then $V$ is $\Ga$-graded by the eigenspaces of $\si$:
$a \in V_\al$ iff $\si a = e^{2\pi\i\al} a$.
For this $\Ga$-grading, Eq.\ \eqref{e4.4} is equivalent to the
``monodromy'' condition $Y(a,\z) = Y(\si a, e^{2\pi\i} \z)$.
Then the above notion of a $\Ga$-twisted module coincides
with the notion of a $\si$-\emph{twisted module} from \cite{FFR,D2}.

\skp
\Remark{ \ref{s4.2}}
Written in terms of modes, Eq.\ \eqref{e4.5} becomes again
Borcherds identity \eqref{e3.28}.
In the definition of a twisted module,
the Borcherds identity (or Jacobi identity) can be replaced 
by associativity relation \eqref{e3.38} (see \cite{Li4}).

\skp

Comparing Eqs.\ \eqref{e3.19} and \eqref{e4.5}, it is obvious that
the notion of a twisted module is a special case of the notion
of a module over a generalized vertex algebra. This can be stated
more precisely as follows (cf.\ \cite{Li2}).

\skp
\Proposition{ \ref{s4.2}}
{\it{
Let\/ $V$ be a\/ $\Ga$-graded
vertex superalgebra, where\/ $\Ga$ is an additive subgroup of\/
$\CC / \ZZ\,$. Introduce the abelian group\/ $Q = (\ZZ / 2\ZZ) \times \Ga$
and a\/ $Q$-grading on\/ $V$ by
\begin{equation*}
a \in V_{(p,\al)} \quad\text{iff}\quad
p(a) = p \,, \;
a \in V_\al \,,
\end{equation*}
where\/ $p(a)$ denotes the parity of\/ $a$.
Let\/ $P=\{\ga_0\}$ be a one-element set with the trivial
action of\/ $Q$, and let
\begin{equation*}
\ett(\hat\al,\hat\be) = (-1)^{p q} \,, \quad 
\De(\hat\al,\hat\be) = \ZZ \,, \quad 
\De(\hat\al,\ga_0) = \al \,, 
\end{equation*}
for\/ $\hat\al=(p,\al), \, \hat\be=(q,\be) \in Q$.
Then a module over\/ $V$ as a generalized vertex algebra is the same
as a\/ $\Ga$-twisted module over\/ $V$ as a vertex superalgebra.
}}

\skp


\section{Generalized vertex algebras 
associated to additive subgroups of vector spaces}\label{s5}

\subsection{Generalized vertex algebras 
associated to vector spaces}\label{s5.1}

Let $\h$ be a vector space (over $\CC$) equipped with a 
symmetric bilinear form $(\cdot | \cdot)$, which is
not necessarily non-degenerate. Inspired by the construction in \cite{FK}
and the construction of lattice vertex algebras \cite{B,FLM,K,LL},
we construct a generalized vertex algebra $V_\h$ as follows.

Let $\hat\h = \h[t,t^{-1}] \oplus \CC K$ be the 
\emph{Heisenberg algebra};
this is a Lie algebra with the bracket
\begin{equation*}
[\al t^m, \be t^n] = m\de_{m,-n}(\al|\be)K \,, \quad [\al t^m,K]=0 \,.
\end{equation*}
The Heisenberg algebra has a representation with $K=1$
on the \emph{Fock space} $F := S(\h[t^{-1}]t^{-1})$ 
(symmetric algebra), so that elements 
of $\h[t^{-1}]t^{-1}$ act by multiplication and $\h[t] 1 = 0$.
(This representation is irreducible iff the bilinear form $(\cdot | \cdot)$
is non-degenerate.)
Let $\CC[\h]$ be the \emph{group algebra} of $\h$ 
considered as an additive group; as a vector space
$\CC[\h]$ has a basis $\{ e^\al \}_{\al\in\h}$, and the
multiplication is $e^\al e^\be := e^{\al+\be}$.

Our generalized vertex algebra $V_\h$ will be
the $\h$-graded vector space
\begin{equation*}
V_\h := F \tt_\CC \CC[\h] = \bigoplus_{\al\in \h} V_\al \,,
\quad V_\al = F \tt e^\al \,.
\end{equation*}
We take $Q=\h$ (as an additive group), and
$\De(\al,\be) = -(\al|\be) + \ZZ$.
Next, we will define certain parafermion fields on $V_\h$.
For every $\al\in\h$, consider the expression
\begin{equation*}
E_\al(\z) := 
\exp\Bigl( \sum_{n>0} \, \frac1{n} \al_{-n} \z^n \Bigr)
\exp\Bigl( \sum_{n<0} \, \frac1{n} \al_{-n} \z^n \Bigr) \,,
\end{equation*}
where $\al_n \in \End F$ denotes the linear operator representing
$\al t^n \in \hat\h$.
It is easy to see that $E_\al(\z)$ is a well-defined field on $F$,
i.e., it is a linear map from $F$ to $F\llbl\z\rrbl$.
We define the \emph{vertex operators}
\begin{equation*}
Y_\al(z) := E_\al(\z) \tt e^\al z^\al 
\colon V_\h \to V_\h \llb\z\rrb \z^\CC \,,
\end{equation*}
where $e^\al$ acts on $\CC[\h]$ by multiplication and 
$\z^\al$ by $\z^\al e^\be := \z^{(\al|\be)} e^\be$.
Then 
\begin{equation*}
Y_\al(\z) v \in V_{\al+\be} \llbl\z\rrbl \z^{(\al|\be)} \,, \qquad
v \in V_\be \,,
\end{equation*}
which means that 
$Y_\al(\z)$ is a parafermion field on $V_\h$ of degree $\al$ 
(see \deref{s3.1}(a)).
We extend the action of the Heisenberg algebra 
from $F$ to $V_\h$ by letting
$\al_n (v \tt e^\be) 
:= (\al_n v) \tt e^\be + \de_{n,0} (\al|\be) v \tt e^\be$.
Then we introduce the \emph{currents}
\begin{equation*}
\al(\z) := \sum_{n\in\ZZ} \, \al_n \, \z^{-n-1} 
\in (\End V_\h) \llb \z, \z^{-1} \rrb \,,
\end{equation*}
which are fields on $V_\h$ of degree $0$.

\skp
\Proposition{ \ref{s5.1}}
{\it{
Introduce the vacuum vector\/ $\vac := 1 \tt e^0 \in V_\h$
and the translation operator\/ $T \in \End V_\h$ determined uniquely by 
\begin{equation*}
T(1 \tt e^\al) := \al t^{-1} \tt e^\al \,, \quad
[T, \al_n] = -n \, \al_{n-1} \,.
\end{equation*}
Let\/ $Q=\h$ $($as an additive group$)$, and let
\begin{equation}\label{e5.8}
\De(\al,\be) = -(\al|\be) + \ZZ \,, \quad
\ett(\al,\be) = e^{\pi\i (\al|\be)} \,. 
\end{equation}
Then the collection of parafermion fields
\begin{equation*}
Y(\al t^{-1} \tt e^0, \z) = \al(\z) \,, \quad
Y(1 \tt e^\al, \z) = Y_\al(z) 
\qquad (\al\in\h)
\end{equation*}
is translation covariant, local and complete,
and hence it defines a unique structure of a generalized vertex algebra
on\/ $V_\h$.
}}

\skp
The \emph{proof} is the same as for lattice vertex algebras,
using \coref{s3.3}
(see \cite[Sect.\ 5.4]{K}), and is in fact 
simpler because of the absence of the $2$-cocycle $\ep$.
Let us recall here only the following crucial formulas,
which imply locality:
\begin{align*}
[\al(\z), \be(\w)] &= (\al|\be) \, \di_\w \de(\z,\w) \,,
\qquad\qquad \al , \be \in\h \,,
\\ 
[\al(\z),Y_\be(\w)] &= (\al|\be) \, Y_\be(\w) \, \de(\z,\w) \,,
\\ 
Y_\al(\z) \, Y_\be(\w) &= 
\io_{\z,\w} (\z-\w)^{(\al|\be)} E_{\al,\be}(\z,\w) 
\tt e^{\al+\be} \z^\al \w^\be
\,,
\end{align*}
where
\begin{equation*}
E_{\al,\be}(\z,\w) = 
\exp\Bigl( \sum_{n>0} 
\frac1{n} \bigl( \al_{-n} \z^n + \be_{-n} \w^n \bigr)
\! \Bigr)
\exp\Bigl( \sum_{n<0} 
\frac1{n} \bigl( \al_{-n} \z^n + \be_{-n} \w^n \bigr)
\! \Bigr)
\end{equation*}
and $\de(\z,\w)$ is the formal delta-function \eqref{e3.17}.\qed

\skp

As an immediate corollary of Propositions \ref{s4.1} and \ref{s5.1},
we obtain a generalized vertex algebra $V_Q= F \tt_\CC \CC[Q]$, 
associated to any additive subgroup $Q \subseteq \h$,
and generated by the collection of parafermion fields
$\{ \al(\z), Y_\be(\z) \}_{\al\in\h , \be\in Q}$.
When the bilinear form $(\cdot | \cdot)$ is non-degenerate, the
algebra $V_Q$ is \emph{simple}
(i.e., it does not have nontrivial proper ideals), because then
the Fock space $F$ is an irreducible $\hat\h$-module.
In the case when $Q$ is a \emph{rational lattice}
(i.e., a free abelian group of finite rank
with a non-degenerate symmetric bilinear form over $\QQ$), 
the generalized vertex algebra $V_Q$
was constructed in \cite{DL,M}.



%

\subsection{Vertex superalgebras associated to integral 
additive subgroups of vector spaces}\label{s5.2}

In this subsection we will continue to use the notation of \seref{s5.1},
and will assume in addition that $\h$ is finite dimensional.
We fix an additive subgroup $Q\subseteq\h$, which is \emph{integral},
i.e., such that $(\al|\be) \in\ZZ$ for all $\al,\be\in Q$.
Let $P$ be the set of all
$\ga\in\h$ such that $(\al|\ga) \in\ZZ$ for all $\al\in Q$. Then $P$
is an additive subgroup of $\h$ containing $Q$,
and the factor group $P/Q$ is isomorphic to a direct product of
a vector space and a finite abelian group.

By Propositions \ref{s4.1} and \ref{s5.1}, we have
a generalized vertex algebra $V_Q = F \tt_\CC \CC[Q]$
and a family of $V_Q$-modules $V_{\ga+Q} = F \tt_\CC e^\ga \CC[Q]$
for every coset $\ga+Q \in P/Q$ (where $\ga\in P\subseteq\h$).
The integrality assumption implies that $\De(\al,\ga) \in \ZZ$
for all $\al\in Q$, $\ga\in P$ (see Eq.\ \eqref{e5.8}).
Thus we can apply \coref{s4.1} to obtain a 
$Q$-graded vertex superalgebra $V_Q^\ep$ and 
a family of $V_Q^\ep$-modules $V_{\ga+Q}^\ep\,$,
for some $2$-cocycle $\ep\colon P \times P \to \CCC$.

Let us spell out this construction in more detail,
utilizing the results of \seref{s3.6}.
First of all, recall that as a vector space
$V_Q^\ep$ coincides with $V_Q$. We define a structure of a vector
superspace on it by letting the parity of $a\in V_\al$ be 
$p(\al) := (\al|\al)$ mod $2\ZZ$; then $\ett(\al,\al) = (-1)^{p(\al)}$.
We choose a $2$-cocycle
$\ep \colon Q \times Q \to \CCC$
such that $\ett^\ep(\al,\be) = (-1)^{p(\al) p(\be)}$.
It corresponds to a central extension \eqref{e4.1} 
with a canonical invariant
(cf.\ Eqs.\ \eqref{e4.11}, \eqref{e5.8})
\begin{equation*}
\oom(\al,\be) := 
(-1)^{ (\al|\al)(\be|\be) + (\al|\be) } \,, \qquad
\al,\be \in Q \,.
\end{equation*}
Since $Q$ is a free abelian group, such $\ep$ can be constructed
explicitly, for example as in \cite[Remark 5.5a]{K}.
It follows from the results of \seref{s3.6} that,
up to isomorphism, $V_Q^\ep$ does not depend on the choice
of a $2$-cocycle $\ep$.

The vertex superalgebra $V_Q^\ep$ is generated by the 
following fields (cf.\ Eq.\ \eqref{e4.9} and \prref{s5.1}):
\begin{align*}
& Y^\ep(\al t^{-1} \tt e^0, \z) = \al(\z) \,, \qquad \al\in\h
& (\text{\textit{currents}}),
\\ 
& Y^\ep(1 \tt e^\al, \z) = Y_\al(z) \, c_\al \,, \qquad \al\in Q
& (\text{\textit{vertex operators}}),
\end{align*}
where $c_\al |_{V_\be} := \ep(\al,\be) \id_{V_\be}$.
Alternatively, we can define $Y^\ep(1 \tt e^\al, \z)$ as $Y_\al(z)$
but replace everywhere the group algebra $\CC[Q]$ with
the \emph{$\ep$-twisted group algebra} $\CC_\ep[Q]$,
which coincides with $\CC[Q]$ as a vector space but has
multiplication defined by Eq.\ \eqref{e4.5e}.

\skp
\Remark{ \ref{s5.2}}
If the bilinear form $(\cdot | \cdot)$ is non-degenerate, 
then $Q$ is an \emph{integral lattice} and
$P$ is a lattice of the same rank as $Q$, called its \emph{dual lattice}.
Then $V_Q^\ep$ coincides with the 
lattice vertex (super)algebra
from \cite{B,FLM,K,LL}. It is a simple vertex (super)algebra, and
the $V_Q^\ep$-modules $V_\Ga^\ep$ $(\Ga\in P/Q)$ are irreducible.
These modules were constructed in \cite{FLM} by a different method
(see also \cite{LL}). It was proved in \cite{D1} that they form
a complete list of irreducible $V_Q^\ep$-modules up to isomorphism.



\section*{Acknowledgments}
The first author wishes to thank 
Nikolay M.\ Nikolov for inspiring discussions and for
collaboration on \cite{BN}. Bakalov was supported in part by an 
FRPD grant from North Carolina State University.
Kac was supported in part by NSF grant DMS-0501395.


\end{document}